\documentclass[a4paper]{amsart}
\usepackage{amsmath}
\usepackage{amsthm}
\usepackage{amssymb}
\usepackage{fullpage}

\newcommand{\A}{\mathbb{A}}

\newcommand{\Z}{\mathbb{Z}}

\newcommand{\fg}{\mathfrak{g}}
\newcommand{\fu}{\mathfrak{u}}

\newcommand{\Coh}{\mathfrak{Coh}}

\newcommand{\cF}{\mathcal{F}}
\newcommand{\cG}{\mathcal{G}}
\newcommand{\cI}{\mathcal{I}}

\newcommand{\cO}{\mathcal{O}}
\newcommand{\cIC}{\mathcal{IC}}

\newcommand{\cg}[1]{\Coh^G(#1)}
\newcommand{\cgl}[2]{\Coh^G(#1)_{\le #2}}
\newcommand{\cgg}[2]{\Coh^G(#1)_{\ge #2}}
\newcommand{\cb}[1]{\Coh^B(#1)}
\newcommand{\cbl}[2]{\Coh^B(#1)_{\le #2}}
\newcommand{\cbg}[2]{\Coh^B(#1)_{\ge #2}}
\newcommand{\Rep}{\mathfrak{Rep}}

\newcommand{\dg}[1]{\mathcal{D}^G(#1)}
\newcommand{\db}[1]{\mathcal{D}^B(#1)}

\newcommand{\hto}{\hookrightarrow}
\newcommand{\la}{\langle}
\newcommand{\ra}{\rangle}

\newcommand{\bru}[1]{X^\circ_{#1}}
\newcommand{\schu}[1]{X_{#1}}
\newcommand{\obru}[1]{Y_{#1}}

\newcommand{\wsame}[1]{\Pi(#1)}
\newcommand{\wopp}[1]{\Theta(#1)}
\newcommand{\wsamel}[1]{\Pi_L(#1)}
\newcommand{\woppl}[1]{\Theta_L(#1)}

\newcommand{\tsame}[1]{\pi(#1)}
\newcommand{\topp}[1]{\theta(#1)}
\newcommand{\tsamel}[1]{\pi_L(#1)}
\newcommand{\toppl}[1]{\theta_L(#1)}

\DeclareMathOperator{\Spec}{Spec}
\DeclareMathOperator{\codim}{codim}
\DeclareMathOperator{\scod}{scod}
\DeclareMathOperator{\Hom}{Hom}
\DeclareMathOperator{\Ext}{Ext}

\DeclareMathOperator{\cRHom}{\mathit{R}\mathcal{H}\mathit{om}}

\newtheorem{thm}{Theorem}[section]
\newtheorem{lem}[thm]{Lemma}
\newtheorem{prop}[thm]{Proposition}

\theoremstyle{definition}

\theoremstyle{remark}

\newcommand{\ssm}{\smallsetminus}

\title{Staggered sheaves on partial flag varieties}
\author{Pramod N.~Achar}
\address{Department of Mathematics\\
  Louisiana State University\\
  Baton Rouge, LA 70803}
\email{pramod@math.lsu.edu}
\thanks{The research of the first author was partially supported by NSF
grant~DMS-0500873.}
\author{Daniel S.~Sage}
\email{sage@math.lsu.edu}
\thanks{The research of the second author was partially supported by NSF
grant~DMS-0606300.}

\begin{document}

\begin{abstract}
\emph{Staggered $t$-structures} are a class of $t$-structures on derived
categories of equivariant coherent sheaves.  In this note, we show that the
derived category of coherent sheaves on a partial flag variety, equivariant
for a Borel subgroup, admits an artinian staggered $t$-structure. As a
consequence, we obtain a basis for its equivariant $K$-theory consisting
of simple staggered sheaves.

%
\end{abstract}

\maketitle

Let $X$ be a variety over an algebraically closed field, and let $G$ be an
algebraic group acting on $X$ with
finitely many orbits. Let $\cg X$ be the category of $G$-equivariant
coherent sheaves on $X$, and let $\dg X$ denote its bounded derived
category. \emph{Staggered sheaves}, introduced in~\cite{stag}, are the objects in the heart of a
certain $t$-structure on $\dg X$, generalizing the perverse coherent $t$-structure~\cite{bez:pc}.  The
definition of this $t$-structure depends on the following data: (1)~an
\emph{$s$-structure} on $X$ (see below); (2)~a choice
of a Serre--Grothendieck dualizing complex $\omega_X \in \dg
X$~\cite{hart}; and (3)~a
\emph{perversity}, which is an integer-valued function on the set of
$G$-orbits, subject to certain constraints.  When the perversity is
``strictly monotone and comonotone,'' the category of staggered sheaves is
particularly nice: every object has finite length, and every simple object
arises by applying an intermediate-extension (``IC'') functor to an
irreducible vector bundle on a $G$-orbit.

An $s$-structure on $X$ is a collection of full subcategories $(\{\cgl
Xn\}, \{\cgg Xn\})_{n \in \Z}$, satisfying various conditions
involving $\Hom$- and $\Ext$-groups, tensor products, and short exact
sequences.  
The \emph{staggered codimension} of the closure of an orbit $i_C: C \to X$,
denoted $\scod \overline C$, is defined to be $\codim \overline C + n$,
where $n$ is the unique integer such that $i_C^!\omega_X \in \dg C$ is a
shift of an object in $\cgl Cn \cap \cgg Cn$.  By~\cite[Theorem~9.9]{stag},
a sufficient condition for the existence of a strictly monotone and
comonotone
perversity is that staggered codimensions of neighboring orbits differ by
at least $2$.  The goal of this note is to establish the existence of a
well-behaved staggered category on partial flag varieties, by constructing
an $s$-structure and computing staggered codimensions.  As a consequence,
we obtain a basis for the equivariant $K$-theory $K^B(G/P)$
consisting of simple staggered sheaves.

\section{A gluing theorem for $s$-structures}
\label{sect:glue}

If $X$ happens to be a single
$G$-orbit, $s$-structures on $X$ can be described via the equivalence
between $\cg X$ and the category of finite-dimensional representations of
the isotropy group of $X$.  In the general case, however, specifying an
$s$-structure on $X$ directly can be quite arduous.  The following
``gluing theorem'' lets us specify an $s$-structure
on $X$ by specifying one on each $G$-orbit.

\begin{thm}\label{thm:sstruc}
For each orbit $C \subset X$, let $\cI_C \subset \cO_X$ denote the ideal sheaf corresponding to the closed subscheme $i_C: \overline C \hto X$.  Suppose each orbit $C$ is endowed with an $s$-structure, and that $i_C^*\cI_C|_C \in \cgl C{-1}$.  There is a unique $s$-structure on $X$ whose restriction to each orbit is the given $s$-structure.
\end{thm}
\begin{proof}
This statement is nearly identical to~\cite[Theorem~10.2]{stag}.  In that result, the requirement that $i_C^*\cI_C|_C \in \cgl C{-1}$ is replaced by the following two assumptions:
\begin{enumerate}
\renewcommand{\labelenumi}{(F\arabic{enumi})}
\item For each orbit $C$, $i_C^*\cI_C|_C \in \cgl C0$.\label{f:ideal}
\item Each $\cF \in \cgl Cw$ admits an extension $\cF_1 \in \cg {\overline
C}$ whose restriction to any smaller orbit $C' \subset \overline C$ is in
$\cgl{C'}w$.\label{f:extend}
\end{enumerate}
Condition~(F1) is trivially implied by the stronger assumption that
$i_C^*\cI_C|_C \in \cgl C{-1}$.  It suffices, then, to show that~(F2) is
implied by it as well.  Given $\cF \in \cgl Cw$, let $\cG \in \cg
{\overline C}$ be some sheaf such that $\cG|_C \simeq \cF$.  Let $C'
\subset \overline C \ssm C$ be a maximal orbit (with respect to the closure
partial order) such that $i_{C'}^*\cG|_{C'} \notin \cgl {C'}w$.  (If
there is no such $C'$, then $\cG$ is the desired extension of $\cF$, and
there is nothing to prove.)  Let $v \in \Z$ be such that $i_{C'}^*\cG|_{C'}
\in \cgl {C'}v$.  By assumption, we have $v > w$.  Let $\cG' = \cG \otimes
\cI_{C'}^{\otimes v-w}$.  Since $\cI_{C'}|_{X \ssm \overline C'}$ is
isomorphic to the structure sheaf of $X \ssm \overline C'$, we see that
$\cG'|_{\overline C \ssm \overline C'} \simeq \cG|_{\overline C \ssm
\overline C'}$.  On the other hand, according to~\cite[Axiom~(S6)]{stag} (which describes how tensor products behave with respect to $s$-structures), the fact that $i_{C'}^*\cI_{C'}|_{C'} \in
\cgl{C'}{-1}$ implies that $i_{C'}^*\cG'|_{C'} \simeq
i_{C'}^*\cG|_{C'} \otimes (i_{C'}^*\cI_{C'}|_{C'})^{\otimes v-w} \in
\cgl{C'}w$.  Thus, $\cG'$ is a new extension of $\cF$ such that the number
of orbits in $\overline C \ssm C$ where~(F2) fails is fewer than for
$\cG$. Since the total number of orbits is finite, this construction can be
repeated until an extension $\cF_1$ satisfying~(F2) is obtained.
\end{proof}

\section{Torus actions on affine spaces}

In this section, we consider coherent sheaves on an affine space.  Let $T$ be an algebraic torus over an algebraically closed field $k$, and let $\Lambda$ be its weight lattice.  Choose a set of weights $\lambda_1, \ldots, \lambda_n \in \Lambda$.  Let $T$ act linearly on $\A^n = \Spec k[x_1,\ldots,x_n]$ by having it act with weight $\lambda_i$ on the line defined by the ideal $(x_j : j \ne i)$.
Given $\mu \in \Lambda$, let $V(\mu)$ denote the one-dimensional
$T$-representation of weight $\mu$.  If $X$ is an affine space with a
$T$-action, we denote by $\cO_X(\mu)$ the twist of the
structure sheaf of $X$ by $\mu$.

Suppose $m \le n$, and identify $\A^m$ with the closed subspace of $\A^n$ defined by the ideal $(x_j : j > m)$.  Let $\cI \subset \cO_{\A^n}$ denote the corresponding ideal sheaf, and let $i: \A^m \hto \A^n$ be the inclusion map.

\begin{prop}\label{prop:torus}
With the above notation, we have
\[
i^*\cI \simeq \cO_{\A^m}(-\lambda_{m+1}) \oplus \cdots \oplus \cO_{\A^m}(-\lambda_n)
\qquad\text{and}\qquad
i^!\cO_{\A^n}(\mu) \simeq \cO_{\A^m}(\mu + \lambda_{m+1} + \cdots
\lambda_n)[m-n].
\]
\end{prop}
\begin{proof}
Throughout, we will pass freely between coherent sheaves and modules, and
between ideal sheaves and ideals.  In the $T$-action on the ring $R =
k[x_1, \ldots, x_n]$, $T$ acts on the one-dimensional space $kx_i$ with
weight $-\lambda_i$.  We have $i^*\cI \simeq \cI/\cI^2 \simeq
(x_{m+1}, \ldots, x_n)/(x_ix_j : m+1 \le i < j \le n)$, so if we let $S =
k[x_1,\ldots, x_m]$, we obtain
$
i^*\cI  
\simeq x_{m+1}S \oplus \cdots \oplus x_n S 
\simeq V(-\lambda_{m+1}) \otimes S \oplus \cdots \oplus
V(-\lambda_n) \otimes S
$.

To calculate $i^!\cO_{\A^n}(\mu)$, we may assume that $m = n-1$, as the
general case then follows by induction.  Recall that $i_*i^!(\cdot) \simeq
\cRHom(i_*\cO_{\A^{n-1}}, \cdot)$.  To compute the latter functor, we
employ the projective resolution $x_nR \hto R$ for $i_*\cO_{\A^{n-1}}$. 
Now, $x_nR \simeq V(-\lambda_n) \otimes R$, so when we apply $\Hom(\cdot,
V(\mu) \otimes R)$ to this sequence, we obtain an injective map
$V(\mu) \otimes R \to V(\mu+\lambda_n) \otimes R$
whose image is $V(\mu+\lambda_n) \otimes x_nR$.  The cohomology of this
complex vanishes except in degree $1$, where we find $V(\mu+\lambda_n)
\otimes R/x_nR$.  Thus, $i_*i^!\cO_{\A^n}(\mu) \simeq
\cRHom(i_*\cO_{\A^{n-1}}, \cO_{\A^n}(\mu)) \simeq
i_*\cO_{\A^{n-1}}(\mu+\lambda_n)[-1]$, as desired.
\end{proof}

\section{$s$-structures on Bruhat cells}

Let $G$ be a reductive algebraic group over an algebraically closed field,
and let $T \subset B \subset P$ be a maximal torus, a Borel subgroup, and a
parabolic subgroup, respectively, and let $L$ be the Levi subgroup of $P$.

Let $W$ be the Weyl group of $G$ (with respect to $T$), and let $\Phi$ be
its root system.  Let $\Phi^+$ be the set of positive roots
corresponding to $B$.  Let $W_L \subset W$ and $\Phi_L \subset \Phi$ be the
Weyl group and root system of $L$, and let $\Phi_P = \Phi_L \cup \Phi^+$. 
For each $w \in W$, we fix once and for all a representative in $G$, also
denoted $w$.
Let
$\bru w$ denote the Bruhat cell $BwP/P$, let $\schu w$ denote its
closure (a Schubert variety), and let $i_w: \schu w \to G/P$ be the
inclusion. Note that $\bru w = \bru v$ if and only if $wW_L = vW_L$.

Let $\Lambda$ denote the weight lattice of $T$, and let $\rho =
\frac{1}{2}\sum \Phi^+$.  (For a set $\Psi \subset \Phi$, we write ``$\sum
\Psi$'' for $\sum_{\alpha \in \Psi} \alpha$.)  For any $w \in W$, we define
various subsets of $\Phi^+$ and elements of $\Lambda$ as follows:
\begin{align*}
\wsame w &= \Phi^+ \cap w(\Phi^+) &
\tsame w &= \textstyle\sum \wsame w & 
&& 
\wsamel w &= \Phi^+ \cap w(\Phi^+ \ssm \Phi_L) &
\tsamel w &= \textstyle\sum \wsamel w \\
\wopp w &= \Phi^+ \cap w(\Phi^-) &
\topp w &= \textstyle\sum \wopp w &
&&
\woppl w &= \Phi^+ \cap w(\Phi^- \ssm \Phi_L) &
\toppl w &= \textstyle\sum \woppl w
\end{align*}

For any subset $\Psi \subset \Phi$, we define $\fg(\Psi) = \bigoplus_{\alpha \in \Psi} \fg_\alpha$.
Next, let $B_w = w Bw^{-1}$, and let $U_w$ denote the unipotent radical of
$B_w$.  Its Lie algebra $\fu_w$ is described by $\fu_w = \fg(w(\Phi^+))$.
Let $\la\cdot,\cdot\ra$ denote the Killing form.
By rescaling if necessary, assume that $\la 2\rho, \lambda \ra \in \Z$
for all $\lambda \in \Lambda$.

Now, the category $\cb{\bru w}$ is equivalent to the category $\Rep(B_w \cap B)$ of representations of the isotropy group $B_w \cap B$.  We define an $s$-structure on $\bru w$ via this equivalence as follows:
\begin{equation}\label{eqn:s-bru}
\begin{aligned}
\cbl{\bru w}n &\simeq \{ V \in \Rep(B_w \cap B) \mid \text{$\la \lambda,
-2\rho \ra \le n$ for all weights $\lambda$ occurring in $V$} \} \\
\cbg{\bru w}n &\simeq \{ V \in \Rep(B_w \cap B) \mid \text{$\la \lambda,
-2\rho \ra \ge n$ for all weights $\lambda$ occurring in $V$} \}
\end{aligned}
\end{equation}

\begin{lem}\label{lem:bru-orb}
For any $v, w \in W$, there is a $B_v$-equivariant isomorphism
$B_vwP/P \simeq \fg(v(\woppl{v^{-1}w}))$.
\end{lem}
\begin{proof}
We have
$
B_v wP/P 
= w \cdot B_{w^{-1}v}P/P \simeq w \cdot B_{w^{-1}v}/(B_{w^{-1}v} \cap P)
$.
Since $B_{w^{-1}v} \cap P$ contains the maximal torus $T$, the quotient $B_{w^{-1}v}/(B_{w^{-1}v} \cap P)$ can be identified with a quotient of $U_{w^{-1}v}$, and hence of $\fu_{w^{-1}v}$.  Specifically, it is isomorphic to $\fg(w^{-1}v(\Phi^+) \ssm \Phi_P) \simeq \fg(w^{-1}v(\Phi^+) \cap (\Phi^- \ssm \Phi_L))$, so
\[
B_vwP/P
\simeq w \cdot \fg(w^{-1}v(\Phi^+) \cap (\Phi^- \ssm \Phi_L)) 
\simeq \fg(v(\woppl{v^{-1}w})).\qedhere
\]
\end{proof}

In the special case $v = ww_0$, where $w_0$ is the longest element of $W$,
the set $v(\woppl{v^{-1}w})$ is given by 
\[
ww_0(\woppl{w_0}) = w(\Phi^-) \cap w(\Phi^- \ssm \Phi_L) = w(\Phi^- \ssm \Phi_L) = -\wsamel w \sqcup \woppl w.
\]
Let $\obru w = B_{ww_0}wP/P$.  Applying Lemma~\ref{lem:bru-orb} with $v = 1$ and with $v = ww_0$, we obtain
\begin{equation}
\bru w \simeq \fg(\woppl w)
\qquad\text{and}\qquad
\obru w \simeq \bru w \oplus \fg(-\wsamel w).
\end{equation}

Finally, let $\cI_w$ denote the ideal sheaf on $G/P$ corresponding to
$\schu w$.  Since $\obru w$ is open, Proposition~\ref{prop:torus} tells us
that $i_w^*\cI_w|_{\bru w} \simeq
\bigoplus_{\alpha \in \wsamel w} \cO_{\bru w}(\alpha)$.  Since $\la \alpha,
-2\rho \ra < 0$ for all $\alpha \in \Phi^+$, we see that
$i_w^*\cI_w|_{\bru w} \in \cbl{\bru w}{-1}$, and then
Theorem~\ref{thm:sstruc} gives us an $s$-structure on $G/P$.  Separately,
Proposition~\ref{prop:torus} also tells us that $i_w^!\cO_{G/P}[\codim
\schu w]$ is in
$\cbl{G/P}{\la \tsamel w, 2\rho\ra} \cap \cbg{G/P}{\la \tsamel w,
2\rho\ra}$.
If $w$ is the unique element of maximal length in its coset
$w W_L$, then we have $\codim \schu w = |\Phi^+| - \ell(w)$ and $\tsamel w
= \tsame w$.  (See~\cite[Chap.~2]{carter}.) Combining these observations
gives us the following theorem.

\begin{thm}\label{thm:scod}
There is a unique $s$-structure on $G/P$ compatible with those on the
various $\bru w$.  If $w$ is the unique element of maximal length in
$wW_L$, then the staggered codimension of $\schu w$, with respect to the
dualizing complex $\cO_{G/P}$, is given by $\scod \schu w = |\Phi^+| -
\ell(w) + \la \tsame w, 2\rho \ra$.\qed
\end{thm}

\section{Main result}

\begin{thm}\label{thm:main}
With respect to the $s$-structure and dualizing complex of Theorem~\ref{thm:scod},
$\db{G/P}$ admits an artinian staggered $t$-structure.  In particular, the
set of simple staggered sheaves $\{\cIC(\schu w, \cO_{\bru w}(\lambda))\}$,
where
$\lambda \in \Lambda$, and $w$ ranges over a set of coset representatives
of $W_L$, forms a basis for $K^B(G/P)$.
\end{thm}

By the remarks in the introduction, this theorem follows from
Proposition~\ref{prop:codim} below.  Throughout this section, the notation
``$u \cdot v$'' for the product of $u, v \in W$ will be used to indicate
that $\ell(uv) = \ell(u) + \ell(v)$.
Note that if $s$ is a simple reflection corresponding to a simple
root $\alpha$, $\ell(sw) > \ell(w)$ if and only if $\alpha \in \wsame w$.

\begin{lem}\label{lem:refl}
Let $s$ be a simple reflection, and let $\alpha$ be the corresponding simple root.  If $\ell(sw) > \ell(w)$, then $\tsame{sw} = s\tsame w + \alpha$ and $\topp{sw} = s\topp w + \alpha$.
\end{lem}
\begin{proof}
Since $\wsame s = \Phi^+ \ssm\{\alpha\}$, it is easy to see that if $\alpha
\in \wsame w$, then $\wsame{sw} = s(\wsame w \ssm \{\alpha\})$, and hence
that $\tsame{sw} = s(\tsame w - \alpha) = s\tsame w + \alpha$.  The proof
of the second formula is similar.
\end{proof}

\begin{lem}\label{lem:prod}
For any $w \in W$, we have $\la \tsame w, \topp w \ra = 0$.
\end{lem}
\begin{proof}
We proceed by induction on $\ell(w)$.  If $w = 1$, $\topp w = 0$, and the
statement is trivial.  If $\ell(w) \ge 1$, write $w = s \cdot v$ with $s$ a
simple reflection.  Let $\alpha$ be the corresponding simple root.  We
have $\la \tsame w, \topp w \ra = \la \tsame{sv}, \topp{sv} \ra
= \la s\tsame v + \alpha, s\topp v + \alpha \ra$, and so
\[
\la \tsame w, \topp w \ra
= \la s\tsame v, s \topp v \ra + \la s \tsame v, \alpha \ra + \la s \topp
v, \alpha \ra + \la \alpha, \alpha \ra
= \la \tsame v, \topp v \ra + \la s(2\rho) + \alpha, \alpha \ra.
\]
Now, $\la \tsame v, \topp v \ra$ vanishes by assumption.  Since $s$
permutes $\Phi^+ \ssm \{\alpha\}$, and $2\rho - \alpha$ is the sum of all
roots in $\Phi^+ \ssm \{\alpha\}$, we see that $s(2\rho - \alpha) = 2\rho
- \alpha$.  But $s(2\rho - \alpha) = s(2\rho) + \alpha$ as well, so we find
that
\[
\la \tsame w, \topp w \ra =
\la 2\rho - \alpha, \alpha\ra = \la s(2\rho - \alpha), \alpha \ra
= \la 2\rho - \alpha, s\alpha\ra = -\la 2\rho - \alpha, \alpha\ra.
\]
Comparing the second and last terms above, we see that all these quantities
vanish, as desired.
\end{proof}

\begin{prop}\label{prop:neg}
If $\alpha \in \wsame w$ is a simple root, then $\la \alpha, \topp w \ra \le 0$.
\end{prop}
\begin{proof}
It is clear that it suffices to consider the case where $W$ is irreducible.  
We proceed by induction on $\ell(w)$.  When $w = 1$, $\topp w = 0$, so the
statement holds trivially.  Now, suppose $\ell(w) > 0$, and let $t$ be a
simple reflection such that $\ell(tw) < \ell(w)$.  Let $\beta$ be the
simple root corresponding to $t$.  We must consider four cases, depending
on the form of $tw$.

{\it Case 1}. $w = t\cdot v$ with $\alpha \in \wsame v$.
Then
$
\la \alpha, \topp{tv} \ra
= \la \alpha, t \topp v + \beta \ra
= \la t \alpha, \topp v \ra + \la \alpha,\beta \ra
$,
so
$
\la \alpha, \topp{tv} \ra
= \la \alpha - \la \beta^\vee,\alpha \ra \beta, \topp v \ra + \la \alpha,\beta \ra 
= \la \alpha, \topp v \ra - \la \beta^\vee, \alpha \ra \la \beta, \topp
v\ra + \la \alpha, \beta \ra
$.
We know that $\la \beta^\vee,\alpha \ra \le 0$ and $\la \alpha,\beta\ra
\le 0$.  The fact that $\ell(tv) > \ell(v)$ implies that $\beta \in \wsame
v$, and $\alpha \in \wsame v$ by assumption, so $\la \alpha, \topp v \ra
\le 0$ and $\la \beta, \topp v\ra \le 0$ by induction.  The result follows.

In the remaining cases, we will have $\alpha \notin \wsame{tw}$.  This implies that $s$ and $t$ do not commute.  Let $N = \la \alpha^\vee,\beta\ra \la\beta^\vee,\alpha \ra$.  We then have $N \in \{1,2,3\}$, with $N = 3$ occurring only in type $G_2$.

{\it Case 2}. $w = ts \cdot v$ with $\beta \in \wsame v$.  We have
$
\la \alpha, \topp{tsv} \ra
= \la \alpha, t \topp {sv} + \beta \ra
= \la \alpha, ts \topp v + t\alpha + \beta \ra
= \la st \alpha, \topp v \ra + \la \alpha, t\alpha + \beta \ra
$.
It is easy to check that $st \alpha = (N-1)\alpha - \la \beta^\vee,\alpha\ra \beta$, and hence that
$
\la st\alpha, \topp v\ra  = (N-1) \la \alpha, \topp v\ra - \la \beta^\vee,
\alpha \ra \la \beta, \topp v\ra
$.
Now, $\beta \in \wsame v$ by assumption, and $\alpha \in \wsame v$ since $\ell(sv) > \ell(v)$, so $\la \alpha, \topp v\ra \le 0$ and $\la \beta, \topp v\ra \le 0$ by induction.  Clearly, $N-1 \ge 0$ and $\la \beta^\vee,\alpha\ra < 0$, so $\la st \alpha, \topp v \ra \le 0$.
Next, we have $t\alpha + \beta = \alpha - \la \beta^\vee, \alpha \ra \beta
+ \beta$, so
$
\textstyle
\la \alpha, t\alpha + \beta \ra = \la \alpha,\alpha\ra - \la \beta^\vee, \alpha\ra\la \alpha, \beta \ra + \la\alpha, \beta\ra
= \frac{\la \alpha,\alpha \ra}{2}(2 -N + \la \alpha^\vee,\beta\ra)
$.
Recall that $\la \alpha^\vee,\beta \ra \in \{-1,-N\}$, so $(2 -N + \la
\alpha^\vee,\beta\ra)$ is either $1-N$ or $2-2N$.  In either case, we see
that $\la \alpha, t\alpha + \beta \ra \le 0$.  It follows that $\la \alpha,
\topp w \ra \le 0$.

In the last two cases, we assume that $\beta \notin \wsame{stw}$.  This implies that $w = tst \cdot v$ for some $v$.  We also have $sw = stst \cdot v$, so it must be that $N \ge 2$.

{\it Case 3}. $w = tst \cdot v$ and $N = 2$.  In this case, $sw = stst \cdot v = tsts \cdot v$, so $\ell(sv) > \ell(v)$, and hence $\alpha \in \wsame v$.  Calculations similar to those above yield that $\topp{tstv} = tst \topp v + ts \beta + t\alpha + \beta$, and that
$
\textstyle
\la \alpha, ts \beta + t\alpha + \beta \ra 
= \la \alpha, \beta \ra 
- \frac{\la \alpha,\alpha\ra}{2}\la \alpha^\vee,\beta\ra = 0
$.
Thus,
$
\la \alpha, \topp{tstv} \ra = \la \alpha, tst \topp v \ra + \la \alpha, ts \beta + t\alpha + \beta \ra = \la tst \alpha, \topp v \ra
$.
Direct calculation shows that $tst \alpha = \alpha$ (regardless of whether $\alpha$ is a short root or a long root).  Since $\alpha \in \wsame v$, $\la \alpha, \topp v\ra \le 0$ by induction, so $\la \alpha, \topp w\ra \le 0$ as well.

{\it Case 4}. $w = tst \cdot v$ and $N = 3$.  Since we have assumed that
$W$ is irreducible, $W$ must be of type $G_2$.  Since $sw = stst \cdot v$,
we must have $v \in \{1,s,st\}$, since $ststst$ is the longest word in $W$.
 First suppose $v = st$.  Since $sw$ is the longest word, we have $\wsame w
= \{ \alpha \}$, and hence $\topp w = 2\rho - \alpha$, so
Lemma~\ref{lem:refl} implies that $\la\alpha, \topp w\ra = 0$.  If $v = s$,
direct calculation gives $\topp w = 2\rho - \alpha - s\beta$, and then
that $\la \alpha, \topp w \ra = \la \alpha, \beta \ra < 0$.  Finally, if $v
= 1$, we find that $\topp w = 2\rho - \alpha - s\beta - st\alpha$, and
again $\la \alpha, \topp w \ra < 0$.
\end{proof}

\begin{prop}\label{prop:order}
Let $s$ be a simple reflection, corresponding to the simple root $\alpha$. 
Let $v, w$ be such that $\ell(vsw) = \ell(v) + 1 + \ell(w)$.  Then $\la
\tsame {vw}, 2\rho \ra - \la \tsame{vsw}, 2\rho \ra = (1 - \la
\alpha^\vee, \topp{v^{-1}}\ra)\la w^{-1}\alpha, 2\rho \ra > 0$.
\end{prop}
\begin{proof}
We proceed by induction on $\ell(v)$.  First, suppose that $v = 1$.  Note
that $\topp{v^{-1}} = 0$.  Since $2\rho = \tsame w + \topp w$,
Lemma~\ref{lem:prod} implies that $\la \tsame w, 2\rho \ra = \la \tsame w,
\tsame w \ra$.  Similarly,
\begin{multline*}
\la \tsame{sw}, 2\rho \ra = \la \tsame{sw}, \tsame{sw} \ra
= \la s\tsame w + \alpha, s \tsame w + \alpha \ra \\
= \la s\tsame w, s\tsame w\ra + 2 \la s\tsame w, \alpha \ra + \la \alpha, \alpha \ra 
= \la \tsame w, \tsame w \ra + 2\la \tsame w, s\alpha \ra + \la 2\rho,
\alpha \ra \\
= \la \tsame w, 2\rho \ra - 2\la \tsame w, \alpha \ra + \la \tsame w +
\topp w, \alpha \ra
= \la \tsame w, 2\rho \ra - \la \tsame w - \topp w, \alpha \ra.
\end{multline*}
It is easy to see that $\tsame w - \topp w = w(2\rho)$, whence it follows
that $\la \tsame w, 2\rho \ra - \la \tsame{sw}, 2\rho \ra = \la
w^{-1}\alpha, 2\rho \ra$.  Finally, the fact that $\ell(sw) > \ell(w)$
implies that $w^{-1}\alpha \in \Phi^+$, so $\la w^{-1}\alpha,2\rho \ra >
0$.

Now, suppose $\ell(v) \ge 1$, and write $v = t \cdot x$, where $t$ is a simple reflection with simple root $\beta$.  Using the special case of the proposition that is already established, we find
\[
\la \tsame{xsw}, 2\rho \ra - \la \tsame{txsw}, 2\rho \ra = \la
w^{-1}sx^{-1}\beta, 2\rho \ra
\qquad\text{and}\qquad
\la \tsame{x w}, 2\rho \ra - \la \tsame{tx w}, 2\rho \ra = \la w^{-1}
x^{-1}\beta, 2\rho \ra.
\]
Combining these with the fact that $s x^{-1}\beta = x^{-1}\beta - \la \alpha^\vee, x^{-1}\beta\ra \alpha$, we find
\begin{multline*}
\la \tsame{txw}, 2\rho \ra - \la \tsame{txsw},2\rho \ra =
(\la \tsame{xw}, 2\rho \ra - \la \tsame{xsw}, 2\rho \ra)
+ ( \la w^{-1}sx^{-1}\beta, 2\rho \ra - \la w^{-1} x^{-1}\beta, 2\rho \ra)
\\
= (1-\la \alpha^\vee, \topp{x^{-1}}\ra) \la w^{-1}\alpha, 2\rho\ra
- \la \alpha^\vee, x^{-1}\beta \ra\la w^{-1}\alpha, 2\rho\ra 
= (1- \la \alpha^\vee, \topp{x^{-1}} + x^{-1}\beta \ra) \la w^{-1}\alpha,
2\rho \ra.
\end{multline*}
An argument similar to that of Lemma~\ref{lem:refl} shows that
$\topp{x^{-1}} + x^{-1}\beta = \topp{x^{-1}t} = \topp{v^{-1}}$, so the
desired formula is established.  Since $\ell(vs) > \ell(v)$, we also have
$\ell(sv^{-1}) > \ell(v^{-1})$, and then Proposition~\ref{prop:neg} tells
us that $\la \alpha^\vee, \topp{v^{-1}} \ra \le 0$.  Thus, $\la \tsame{vw},
2\rho \ra - \la \tsame{vsw},2\rho \ra > 0$.
\end{proof}

The preceding proposition is a statement about a pair of adjacent elements
with respect to the Bruhat order.  It immediately implies that for any $v,
w \in W$ with $v < w$ in the Bruhat order, $\la \topp v, 2\rho \ra - \la
\topp w, 2\rho \ra > 0$.  By Theorem~\ref{thm:scod}, we deduce the
following result, and thus establish Theorem~\ref{thm:main}.

\begin{prop}\label{prop:codim}
If $\schu v \subset \overline{\schu w}$, then $\scod \schu v - \scod \schu w \ge 2$. \qed
\end{prop}

\end{document}